 \documentclass[leqno]{amsart}
  \usepackage{amsfonts, amsmath}
 \usepackage{mathrsfs}
\usepackage[colorlinks=false,linktocpage=true]{hyperref}
\usepackage{tocvsec2}
\usepackage{latexsym}
\def\qed{\quad\vcenter{\hrule\hbox{\vrule height .6em\kern .6em \vrule}\hrule}}
\mathchardef\Tau='124
\mathchardef\Iota='111
\newtheorem{thm}{Theorem}[section]

\newtheorem{lem}{Lemma}[section]

\newtheorem{rem}{Remark}[section]
\newtheorem{notn}{Notation}
\newcommand{\thmref}[1]{Theorem~\ref{#1}}

\newcommand{\remref}[1]{Remark~\ref{#1}}
\def\qed{\quad\vcenter{\hrule\hbox{\vrule height.6em\kern.6em\vrule}\hrule}}

\newenvironment{pf*}[1]{{\bigskip\textit{\newline#1.}\quad}}{$\qed$\bigskip\newline}
\numberwithin{equation}{section}
\def\dstyle{\displaystyle}

\def\eqdef{\stackrel{\triangle}{=}}
\def\utx{u(t,x)}
\def\ux{u(x)}
\def\sT{{\mathscr T}^x}

\def\TBTG{\Tau_G^x}

\def\Xx{X^x}

\def\BM{B(t)}
\def\RBM{|B(t)|}

\def\BT{{\Bbb X}^x_{B}}
\def\BTP{{\Bbb X}^x_{B}(t)}

\def\EBT{{\Bbb X}^x_{B,e}}
\def\EBTP{{\Bbb X}^x_{B,e}(t)}
\def\kEBT{{\Bbb X}^{x,k}_{B,e}}
\def\kEBTP{{\Bbb X}^{x,k}_{B,e}(t)}
\def\oEBTP{{\Bbb X}^{x,1}_{B,e}(t)}
\def\BTPstp{{\Bbb X}^x_{B}(\TBTG)}

\def\P{{\Bbb P}}

\def\EP{{\Bbb E}_{\P}}

\def\N{{\Bbb N}}
\def\R{{\Bbb R}}
\def\Rd{{\Bbb R}^d}
\def\pR{{\Bbb R}_+}
\def\OFP{(\Omega,{\mathscr F},\{\Filt\},\P)}
\def\Filt{{\mathscr F}_t}
\def\SG{\mathscr T}
\def\GN{\mathscr A}
\def\GNRt{{\mathscr A}^*_t}

\def\GNRy{{\mathscr A}^*_y}

\def\HDGNs{\mathscr A^{1/2}_s}
\def\eqlaw{\overset{\mathscr L}{=}}

\def\eil1{e_i^\clubsuit(1)}

\def\el1{e^\clubsuit_1}
\begin{document}
  \title[Brownian-Time Processes and PDEs]{Brownian-Time Processes: The PDE Connection and the Half-Derivative Generator}%
  \author{Hassan Allouba and Weian Zheng}
  \address{Department of Mathematics \\
University of Massachusetts Amherst then Indiana University \\
Rawles Hall \\
831 E. 3rd Street \\
Bloomington, IN 47405-7106 \\
E-mail: allouba@indiana.edu}
\address{Department of Statistics\\
East China Normal University, Shanghai,\\ China
\\ and \\
Department of Mathematics \\
 University of California Irvine \\
Irvine, California 92697-3875 \\
E-mail: wzheng@math.uci.edu}

  \date{October 1999 (Original version)  January 2001 (This version)}
    \maketitle
  \begin{abstract}%  ABSTRACT
 We introduce a class of interesting stochastic processes based on Brownian-time processes.
These are obtained by taking Markov processes and replacing the time parameter with the
modulus of Brownian motion.
They generalize the iterated Brownian motion (IBM) of Burdzy and the Markov snake
of Le Gall, and they introduce new interesting examples.
After defining Brownian-time processes, we relate them to
fourth order parabolic PDEs.   We then study their exit problem
as they exit nice domains in $\Rd$, and connect it to elliptic PDEs.  We show that
these processes have the peculiar property that they
solve fourth order parabolic PDEs, but their exit distribution---at least in the standard Brownian-time process case---solves the usual
second order Dirichlet problem.  We recover
fourth order PDEs in the elliptic setting by encoding the iterative nature of the
Brownian-time process, through its exit time, in a standard Brownian motion.
We also show that it is possible to assign a formal generator to these non-Markovian
processes by giving such a generator in the half-derivative sense.%
\end{abstract}
   \setcounter{section}{-1}
   \section{Introduction}
Let $B(t)$ be a one-dimensional Brownian motion starting at $0$ and $X^x(t)$
be an independent $\Rd$-valued continuous Markov process started at $x$, both defined on
a probability space $\OFP$.  We call the process
$\BTP\eqdef X^x(|B(t)|)$
a Brownian-time process (BTP).  In the special case where $X^x$ is a Brownian
motion starting at $x$ we call the process $\BTP$  a Brownian-time Brownian
motion (BTBM).  Excursions-based Brownian-time processes (EBTPs) are obtained
from BTPs by breaking up the path of $\RBM$ into excursion
intervals---maximal intervals $(r,s)$ of time on which $\RBM>0$---and, on each
such interval, we pick
an independent copy of the Markov process $\Xx$ from a finite or an infinite
collection.  BTPs and their close cousins EBTPs may be regarded as canonical
constructions for several famous as well as
interesting new processes.  To see this, observe that the following processes have
the one dimensional distribution $\P(\BTP\in dy)$:
\begin{enumerate}
\renewcommand{\labelenumi}{(\alph{enumi})}
\item Markov snake---when $|B_t|$ increases we generate a new
independent path.  See Le Gall (\cite{LeGall1}, \cite{LeGall2}, and \cite{LeGall3})
for applications to the nonlinear PDE $\Delta u=u^2$.
\item Let $X^{x,1}(t), \ldots,X^{x,k}(t)$ be independent copies of
$\Xx(t)$ starting from point $x$.  On each excursion interval of $\RBM$ use one of the
$k$ copies chosen at random.  When $x=0$, $X^x$ is a Brownian motion starting
at $0$, and $k=2$ this reduces to the iterated Brownian motion (IBM).  See
Burdzy (\cite{Burdzy1} and \cite{Burdzy2}).  We identify such a process by
the abbreviation $k$EBTP and we denote it by $\kEBTP$.  Of course, when $k=1$ we obtain
a BTP.
\item Use an independent copy of $X^x$ on each excursion interval of $\RBM$.
This is the $k\to\infty$ limit of (b) (for a rigorous statement and proof, see the
Appendix).  It is
intermediate between IBM and the Markov snake.
Here, we go forward on a new independent path only after $|B_t|$
reaches $0$.  This process is abbreviated as EBTP and is denoted by $\EBTP$.
\end{enumerate}
In Sections 1 and 2 we connect $\BTP$, $\kEBTP$, and $\EBTP$ to new fourth order parabolic
PDEs and to second and fourth order elliptic PDEs.
As a special case of our results, we get the missing
connection of the IBM of Burdzy to PDEs.  There are of course other iterated
processes that have been linked to different PDEs (see \cite{F}, \cite{H}, and the references
therein), but none solves the IBM PDE.
In Section 3 we show that, eventhough $\BTP$
is not Markovian, we can still assign to it a ``generator'' in the half-derivative sense,
which we therefore call the half-derivative generator.

In Section 1 the PDE connection is given by
\begin{thm}
Let $\SG_sf(x)=\EP f(\Xx(s))$ be the semigroup of the continuous Markov process $\Xx(t)$ and
$\GN$ its generator.  Let $f$ be a bounded measurable function in the domain of
$\GN$, %${\Bbb D}(\GN)$.
with $D_{ij}f$ bounded and H\"older continuous for $1\le i,j\le d$.
If $u(t,x)=\EP f(\kEBTP)$ for any $k\in\N$ $($as stated before $\oEBTP=\BTP$\/$)$,
or if $u(t,x)=\EP f(\EBTP)$,
then $u$ solves the PDE $u$ solves the
\begin{equation}
\begin{cases}
\dstyle \frac{\partial}{\partial t} u(t,x)=\frac{\GN f(x)}{\sqrt{2\pi t}} +
   \frac12  \GN^2u(t,x);
 &t>0,\,x\in\Rd,
\cr u(0,x)=f(x); &x\in\Rd,
\end{cases}
\label{BTPPDE}
\end{equation}
where the operator $\GN$ acts on $u(t,x)$ as a function of $x$ with $t$ fixed.
In particular, if $\BTP$ is a BTBM and
$\Delta$ is the standard Laplacian, then $u$ solves
\begin{equation}
\begin{cases}
  \dfrac{\partial}{\partial t}u(t,x)=\dfrac{\Delta f(x)}{\sqrt{8\pi t}}+
\dfrac18\Delta^2u(t,x);& t>0,\,x\in\Rd,
\cr u(0,x)=f(x);& x\in\Rd.
\end{cases}
\label{BTBMPDE}
\end{equation}
\label{BTPPDES}
\end{thm}
\begin{rem}
The inclusion of the initial function $f(x)$ in the PDEs  \eqref{BTPPDE}
and \eqref{BTBMPDE} is a reflection of the non-Markovian property of our BTP.
Thus the role of $f$ here is fundamentally different from its role in the
standard Markov-PDE connection.
\end{rem}

In Section 2, we focus on BTBMs and we take up the exit problem for $\BTP$.  Towards this end, let $G$ be a bounded open subset of $\Rd$ with regular
boundary $\partial G$.  Each time, we start $\BTP$ at a point $x\in G\cup\partial G$, and we let $\TBTG:=\inf\{t\ge0;\BTP\not\in G\}$.
Our first result says that if we look at the exit distribution of our iterated process we solve the usual {\it second} order Dirichlet problem.  This might seem
surprising at first, but upon reflection, we see that the iterated nature of both our process $\BTP$ and its exit time $\TBTG$ ``cancel'' each other, and we are
effectively reduced to the exit distribution of an ordinary Brownian motion.  For a precise explanation of this phenomenon see
the proof of \thmref{BTPSTP} .  In \thmref{felliptic} we show how to ``recover''
the fourth order PDE in this elliptic setting.
\begin{thm} Let $G$ and $\TBTG$ be defined as above.
If $\ux=\EP f(\BTPstp)$ then $u$ satisfies the Dirichlet problem
\begin{equation}
\begin{cases}
\Delta\ux=0;&x\in G,
\cr u(x)=f(x);&x\in\partial G.
\end{cases}
\label{BTPHAR}
\end{equation}
\label{BTPSTP}
\end{thm}
The next result links the first exit time $\TBTG$ of the BTP $\BTP$ to fourth order PDEs
\begin{thm} Let $G$ and $\TBTG$ be defined as above.
If $\ux=\EP \TBTG$, then $u$ satisfies
\begin{equation}
\begin{cases}
\Delta^2u(x)=8; \qquad &x\in G,
\cr u(x)=0;&x\in\partial G.
\end{cases}
\label{GetoorM}
\end{equation}
\label{EXTTIM}
\end{thm}
We now show how to ``encode'' the
iterated nature of our BTBM process in a Brownian motion so as to recover a fourth order elliptic PDE.
The idea is to look at the Brownian motion $\Xx$ evaluated at the iterated exit time $\TBTG$ (the first
exit time for the iterated process $\Xx(\RBM)$); i.e., $\Xx(\TBTG)$.  Note that
this is {\it not} the exit distribution of $\Xx$ (since $\TBTG\neq\tau_G^x
=\inf\{t\ge0;\Xx\notin G\}$ in general).
The fact that $\TBTG$ is not
a stopping time with respect to the natural filtration of $\Xx$ makes it inconvenient
to deal with directly, so we are led to the deterministic time that captures the desired properties of
$\TBTG$, namely $\sT=\EP\TBTG$.
\begin{rem}
If $x\in\partial G$ then $\TBTG=0$ a.s. $\P$ and so $\sT=\EP\TBTG=0$.
Of course, by \thmref{EXTTIM}, $\sT$ satisfies \eqref{GetoorM}.
We are now ready to give the elliptic fourth order PDE connection to a
Brownian motion at the expected value of the iterated exit time $\TBTG$.
\label{shadow}
\end{rem}
\begin{thm}
Assume that $\Xx$ is the outer Brownian motion in $\BTP=\Xx(\RBM)$, starting at $x$
under $\P$, and let $\sT=\EP\TBTG$.  Let
$f\in C^4(\R^d;\R)$ be biharmonic $($$\Delta^2f\equiv0$\/$);$ and assume polynomial
growth for $f$ and all of its partial derivatives of order $k\le4$.  Then $\ux=\EP f(\Xx(\sT))$ satisfies
\begin{equation}
\begin{cases}
\Delta^2u(x)=4\Delta f(x)+\alpha(x)+\beta(x); \qquad &x\in G,
\cr u(x)=f(x);&x\in\partial G,
\end{cases}
\label{fell2}
\end{equation}
\label{felliptic}
where $\alpha(x)=\nabla(\Delta f(x))\cdot\nabla(\Delta\EP[\tau_G^x]^2)$ and
$\beta(x)=2{\dstyle\sum _{1\leq i,j\leq d\atop i\ne j}D_{ij}
\Delta f(x)D_{ij}\EP[\tau_G^x]^2}$.  In particular, if in addition to the above assumptions on
$f$ and its partial derivatives we assume that $\nabla(\Delta f(x))=0$, where $\nabla$ is the usual gradient,
then $\ux=\EP f(\Xx(\sT))$ solves
\begin{equation}
\begin{cases}
\Delta^2u(x)=4\Delta f(x); \qquad &x\in G,
\cr u(x)=f(x);&x\in\partial G.
\end{cases}
\label{fell1}
\end{equation}
\end{thm}
\begin{rem}
Comparing equation \eqref{fell2} and \eqref{fell1} with \eqref{BTBMPDE},
we see that they all include the bi-Laplacian of $u$ and the Laplacian of the function $f$.
So that, also in the elliptic case \eqref{fell2}, $f$ plays s fundamentally
different role than in the usual Brownian motion-PDE connection: it acts on
$G\cup\partial G$, and not just on the boundary $\partial G$.
\label{rell}
\end{rem}
The following result attaches a formal generator to our BTPs, in the half-derivative sense.  More
precisely, we have
\begin{thm}
Let $X^{x}$ be the outer Markov process in our BTP, starting at $x\in\Rd$ under $\P$.
Suppose that the generator $\GN$ of $X^{x}$ is given by a divergence
form second order partial differential operator as in \eqref{DF}.  Let $\GNRt$ be the generator of
the time-reversed Markov process $\{X^{x}(T-t);0\le t\le T\}$ and suppose that
$C_0^2(\Rd;\R)\supset{\Bbb D}(\GN)\cap{\Bbb D}(\GNRt)$, where ${\Bbb D}(\GN)$ and
${\Bbb D}(\GNRt)$ are the domains of $\GN$ and $\GNRt$, respectively.  Finally, assume that
condition \eqref{zz} holds.  If
\begin{equation}
\HDGNs f(x)\eqdef
\lim_{t\searrow s}\dfrac{\EP[f(\BTP)|\BT(s)]-f(\BT(s))}{{(t-s)}^{1/2}};\ 0<s\le t,
\label{HDGN1}
\end{equation}
then $\HDGNs f(x)$ is given by
\begin{equation}
\begin{split}
{1\over \sqrt{2 \pi }}\left[\GN f(\BT(s))+
\frac{\dstyle\int_0^\infty p(0,s;0,y)h(0,y;x,\BT(s))\GNRy f(\BT(s))dy}
{\dstyle\int_0^\infty p(0,s;0,y)h(0,y;x,\BT(s))dy}\right],
\end{split}
\label{HDGN2}
\end{equation}
where $p(s,t;x,y)$ and $h(s,t;x,y)$ are the transition densities
$($with respect to Lebesgue measure\/$)$ of $|B(t)|$ and $X(t)$, respectively.  In
particular, if $\GN=\GNRt$ for all $t$, then $\HDGNs f(x)$ is
simply ${\sqrt{\dstyle\frac{2}{\pi}}}\GN f(\BT(s)).$
\label{HDGNthm}
\end{thm}
\begin{notn}
We alternate freely between the notations $X(t)$ and $X_t$ for aesthetic reasons
and for typesetting convenience.
\end{notn}
\section{Proof of \thmref{BTPPDES}}
We first prove the theorem for the case of $\utx=\EP f(\BTP)$ using the following
generator computation:
\begin{equation}
\EP f(\BTP) = 2 \int_0^\infty p_t(0,s) \SG_s f(x)ds,
\label{expc}
\end{equation}
where $p_t(0,s)$ is the transition density of $B(t)$.
Differentiating $\eqref{expc}$ with respect to $t$ and putting the derivative under the integral,
which is easily justified by the dominated convergence theorem, then using the
fact that $p_t(0,s)$ satisfies the heat equation we have
\begin{align*}
\begin{split}
\frac{\partial}{\partial t} \EP f(\BTP) & = 2\int_0^\infty\frac{\partial}{\partial t}
p_t(0,s) \SG_s f(x) ds\\
& = \int_0^\infty\frac{\partial^2}{\partial s^2} p_t(0,s) \SG_s f(x) ds
\end{split}
\end{align*}
We now integrate by parts twice, and observe that the boundary terms always vanish
at $\infty$ (as $s\nearrow\infty$) and we have
$\dstyle{(\partial/\partial s) p_t (0,s) =0}$ at $s=0$ but $p_t(0,0) >0$.  Thus,
\begin{align*}
\begin{split}
\frac{\partial}{\partial t} \EP f(\BTP)& = - \int_0^\infty\frac{\partial}{\partial s} p_t (0,s) \frac{\partial}{
\partial s} \SG_s f(x) ds \\
& = p_t (0,0) \GN f(x) + \int_0^\infty p_t (0,s) \GN^2\SG_s f(x) ds
\end{split}
\end{align*}
Taking the application of $\GN^2$ outside the integral and writing
$u(t,x)=\EP f(\BTP)$ we have
$$\frac{\partial}{\partial t} u(t,x) = p_t (0,0) \GN f(x) + \frac12 \GN^2 u(t,x),$$
where, clearly, the operator $\GN$ acts on $u(t,x)$ as a function of $x$ with $t$ fixed.
Obviously, $u(0,x)=f(x)$, so that $\utx=\EP f(\BTP)$ solves \eqref{BTPPDE}.

To prove the result for $\kEBTP$ for $k\in\N\setminus \{1\}$, we show that
$\EP f(\kEBTP)=\EP f(\BTP)$.  Towards this end, let $e^-(t)$ be the
$\RBM$-excursion immediately preceding the excursion straddling $t$, $e(t)$; and
condition on the event that we pick the $j$-th copy of $\Xx$ on $e^-(t)$
(uniformly from among the $k$ available independent copies of $X^x$), using
the independence of the choice of the process $X^{x,j}$ on $e^-(t)$ from $\BM$ and
from the following choice of the $X^x$ copy, on $e(t)$, to get
\begin{equation*}
\begin{split}
\EP f(\kEBTP) &= 2\sum _{j=1}^{k} \int_0^\infty p_t(0,s) \SG_s f(x)
\P[\mbox{we pick the $j$-th copy on } e^-(t)]ds\\
&=\frac2k\sum _{j=1}^{k} \int_0^\infty p_t(0,s) \SG_s f(x)ds
=2\int_0^\infty p_t(0,s) \SG_s f(x)ds\\&=\EP f(\BTP).
\end{split}
\end{equation*}

Finally, to prove that $u(t,x)=\EP f(\EBTP)$ solves \eqref{BTPPDE}, we use the
fact (proven in the Appendix) that $\kEBT\Longrightarrow\EBT$, for some
subsequence $\kEBT$.  Following Skorohod's celebrated result, we may construct
processes $Y_k\eqlaw\kEBT$ and $Y\eqlaw\EBT$ on some probability space such that
$Y_k\longrightarrow Y$ as $k\to\infty$ a.s. uniformly in $t$ on compact sets
of $\pR$.
The result then follows since $\EP f(\kEBTP)=\EP f(\BTP)$ for each $k$ and
since $f$ is bounded and continuous.
$\qed$
\section{Exit PDEs for $\BTP$}
Throughout this section the outer process $\Xx$ is always assumed to be a Brownian motion starting
at $x$ under $\P$, and $G$ is a bounded open subset of $\Rd$ with regular boundary $\partial G$.
  \begin{pf*}{Proof of \thmref{BTPSTP}}
Let
$$\tau_G^{x}\eqdef\inf\{t\ge0;\Xx(t)\notin G\} \mbox { and }\sigma_B^x\eqdef\inf\left\{t\ge0;|B(t)|=\tau_G^{x}\right\},$$
of course $\sigma_B^x=\TBTG$.
We then have
\begin{equation}
\begin{split}
\ux&=\EP f\left[\BTPstp\right]
=\EP f\left[\Xx(\tau_G^{x})\left||B(\sigma_B^x)|=\tau_G^{x}\right.\right]\P\left[|B(\sigma_B^x)|=\tau_G^{x}\right]\\
&=\EP f\left[\Xx(\tau_G^{x})\right],
\end{split}
\label{kexcexp}
\end{equation}
where the last equality in equation \eqref{kexcexp} follows from the obvious fact that $$\P\left[|B(\sigma_B^x)|=\tau_G^{x}\right]=1,$$
a fact which also clearly gives us the independence of the event $\{|B(\sigma_B^k)|=\tau_G^{x}\}$ from $\Xx(\tau_G^{x})$.

Now, $u(x)\eqdef\EP f\left[\Xx(\tau_G^{x})\right]$ is a harmonic function in $G$ (since $\Xx$ is a Brownian motion starting at $x$ under $\P$, and
$\tau_G^{x}$ is its first exit time from G).  It follows that $\ux$ solves the Dirichlet problem \eqref{BTPHAR}.  $\qed$
 \end{pf*}
We then prove the connection of the iterated exit time $\TBTG$ to fourth order PDEs.
  \begin{pf*}{Proof of \thmref{EXTTIM}}
Let $u(x)=\EP\TBTG$ and observe that
\begin{equation}
\begin{split}
\TBTG\eqdef\inf\{t\ge0;\BTP\notin G\}&=\inf\{t\ge0;\RBM\notin[0,\tau_G^x)\}\\
&=\inf\{t\ge0;B(t)\notin(-\tau_G^x,\tau_G^x)\},
\end{split}
\label{equivdef}
\end{equation}
where $\tau_G^x\eqdef\inf\{t\ge0;\Xx(t)\notin G\}$.  Thus, conditioning on
$\tau_G^x$ we easily get
\begin{equation}
u(x)=\EP[\EP[\TBTG|\tau_G^x]]=\EP(\tau_G^x)^2.
\label{Treltau}
\end{equation}
But, from \cite{G} and \cite{Kim} we have that
$u(x)=\EP(\tau_G^x)^2$ solves the equation
\begin{equation}
\Delta^2u=8,
\label{Get}
\end{equation}
for any smooth bounded domain $G$. Plainly,
$u(x)=0$ for $x\in\partial G$.  We thus obtain \eqref{GetoorM}
and this completes the proof.   $\qed$
\end{pf*}
We are now ready to prove \thmref{felliptic}.
\begin{pf*}{Proof of \thmref{felliptic}}
Let
$u(x)=\EP f(\Xx_{\sT})$,
and let $\tau_G^x$ be the first exit time for the Brownian motion
$\Xx$.
It\^{o}'s formula, applied twice gives us
\begin{equation*}
\begin{split}
f(\Xx_{\sT})-f(x)&=\intop _{0}^{\sT}
\nabla f(\Xx_s)\cdot d\Xx_s
+\frac12\intop_{0}^{\sT} \Delta f(\Xx_s)ds\\
&=\intop_{0}^{\sT} \nabla f(\Xx_s)\cdot d\Xx_s
+\frac12 \sT\Delta f(x)+\frac12\intop_{0}^{\sT}\intop_{0}^{s}\nabla(\Delta f(\Xx_r))\cdot
d\Xx_rds\\&+
\frac14\intop_{0}^{\sT}\left[\intop_{0}^{s}\Delta^2f(\Xx_r)dr
\right]ds\\
&=\intop_{0}^{\sT} \nabla f(\Xx_s)\cdot d\Xx_s
+\frac12 \sT\Delta f(x)+\frac12\intop_{0}^{\sT}\intop _{0}^{s}\nabla(\Delta f(\Xx_r))\cdot
d\Xx_rds\\&=\intop_{0}^{\sT} \nabla f(\Xx_s)\cdot d\Xx_s
+\frac12 \sT\Delta f(x)\\&+
\frac12\intop_0^{\sT}\intop_0^{\sT} 1_{\{r<s\}}(r)\nabla (\Delta f(X_r^x))dX_rds\\
&=\intop_{0}^{\sT} \nabla f(\Xx_s)\cdot d\Xx_s
+\frac12 \sT\Delta f(x)+\frac12\intop_0^{\sT}(\sT -r)\nabla (\Delta f(X_r^x))dX_r,
\end{split}
\end{equation*}
where we used the assumption that $\Delta^2f\equiv0$ to get the third equality.
Now, Taking expectations, we get
that all the expectations involving stochastic integrals vanish.  This is because
we assumed that both $\nabla f(x)$ and $\nabla (\Delta f(x))$ have polynomial
growth while the density of $X_r^x$ has exponential decay, so that
\begin{equation*}
\EP\left[\int_0^{\sT}|\nabla f(X^x_s)|^2ds\right]<\infty,\mbox{ and  }
\EP\left[\int_0^{\sT}|(\sT -r)\nabla(\Delta f(X_r^x))|^2dr\right]<\infty.
\end{equation*}
We then have
\begin{equation}
\EP f(\Xx_{\sT})-f(x)=\frac12\sT\Delta f(x).
\label{lexpito}
\end{equation}
Applying the bi-Laplacian to both sides of \eqref{lexpito}; and remembering that
$u(x)=\EP f(\Xx_{\sT})$, that $\Delta^2f\equiv0$, and that $\sT=\EP\TBTG$
(by assumption) and invoking \eqref{Treltau} and \eqref{Get},
we obtain
\begin{equation}
\begin{split}
\Delta^2u(x)&=\frac12\Delta^2[\sT\Delta f(x)]=
\frac12\Delta^2[\sT]\Delta f(x)\\&+\nabla(\Delta f(x))\cdot\nabla(\Delta[\sT])
+\Delta(\nabla(\Delta f(x))\cdot\nabla[\sT])\\
&=4\Delta f(x)+\nabla(\Delta f(x))\cdot\nabla(\Delta[\sT])\\
&+\Delta(\nabla(\Delta f(x))\cdot\nabla[\sT])\\
&=4\Delta f(x)+\nabla(\Delta f(x))\cdot\nabla(\Delta[\sT])\\
&+2\sum _{1\leq i,j\leq d\atop i\ne j}D_{ij}\Delta f(x)D_{ij}[\sT]\\
&=4\Delta f(x)+\nabla(\Delta f(x))\cdot\nabla(\Delta\EP[\tau_G^x]^2)\\
&+2{\dstyle\sum _{1\leq i,j\leq d\atop i\ne j}D_{ij}
\Delta f(x)D_{ij}\EP[\tau_G^x]^2};\quad x\in G,
\end{split}
\label{mypde}
\end{equation}
with the convention that $\sum _{i\ne j}D_{ij}\Delta f(x)D_{ij}\EP[\tau_G^x]^2=0$
if $d=1$.

Finally, as stated in \remref{shadow}, $\sT=0$ whenever $x\in\partial G$, and
so $u(x)=\EP f(\Xx(\sT))=f(x)$ for every $x\in\partial G$.  $\qed$
\end{pf*}
\section{The Half-Derivative Formal generator.}
In this section, we prove the formula for the half-derivative generator of our Brownian-time processes.
We denote by $p(s,t;x,y)$ and $h(s,t;x,y)$ the transition densities
(with respect to Lebesgue measure) of $|B(t)|$ and $X(t)$, respectively.  We denote
the generator of $X$ by $\GN$, and we assume that $X(0)=x_0$ is deterministic.

It is well-known that, for each fixed but arbitrary $0<T<\infty$, the time reversed process
$X^*_T=\left\{X^*_T(t)\eqdef X(T-t);0\leq t\leq T\right\}$ is still
Markovian; we denote its (time-dependent) generator by $\GNRt$.
We assume for simplicity
that $\GN$ is given by a divergence
form second order partial differential operator
\begin{equation}
\GN f=\sum_{i,j=1}^d{\partial\over \partial x_i}
\left[g^{ij}(x){\partial\over \partial x_j}f\right],
\label{DF}
\end{equation}
where $d$ is the space dimension and $g^{ij}\in C^2(\Rd;\R)$ satisfies $c<g^{ij}(x)<c^{-1}$
for some positive constant $c$.  From Aronson's inequality we have a constant
$c_1$ such that
\begin{eqnarray}\label{ar} h(s,t;x,y )\leq {c_1\over (t-s)^{d/2}}
\exp\left\{ -{|x-y|^2\over c_1(t-s)}\right\}.
\end{eqnarray}
Moreover (see, for example, \cite{lz} and \cite{z})
\begin{equation}
\GNRt f=\GN f+2\sum_{i,j=1}^d{\partial\over \partial x_i}\log h(0,t;x_0,x) g^{ij}(x)
{\partial\over \partial x_j}f
\label{GNR}
\end{equation}
In particular, when $\GN=\dstyle{1\over 2}\Delta$,
$\GNRt=\displaystyle{{1\over 2}\Delta +\frac{x_0-x}{ t}}\nabla$.

We assume that for every $f\in C_0^2(\Rd;\R)$
\begin{eqnarray}\label{ass1}
\lim_{t\searrow s}|t-s|^{-1}\left[\int h(s,t;x,y)f(y)dy-f(x)\right]=\GN f(x),
\end{eqnarray}
\begin{equation}\label{ini1}
\lim_{s\nearrow t}|t-s|^{-1}\left[\int{ h(0,s;x_0,y)h(s,t;y,x)f(y)\over h(0,t;x_0,x)}dy-f(x)\right]
=\GNRt f(x),
\end{equation}
and without losing generality we assume that there is a constant $0<c_2<\infty $ such that
\begin{eqnarray}\label{zz}
{\partial\over \partial x_i}\log h(0,t;x_0,x)\leq c_2{|x_0-x|+c_2\over t^{c_2}}.
\end{eqnarray}
When $\GN$ is the Laplacian, the above condition is easily satisfied.
It is easy to deduce
\begin{lem}\label{hhh}
For any fixed $f\in C^2_0(\Rd;\R)$ and $x\in\Rd$, there is a constant
$0<c_3<\infty $ such that
\begin{equation}\label{ini}
\sup_{s<t}\left\{|t-s|^{-1}\left[\int{h(0,s;x_0,y)h(s,t;y,x)f(y)\over
h(0,t;x_0,x)}dy-f(x)\right]\right\}< c_3t^{-c_2}.
\end{equation}
\end{lem}
\begin{pf*}{Proof}
Since $\int h(0,s;x_0,y)h(s,t;y,x)dy= h(0,t;x_0,x)$, then
$$\int{h(0,s;x_0,y)h(s,t;y,x)f(y)\over h(0,t;x_0,x)}dy$$
is bounded by the same bound on $f$.  Thus,
when $s<2^{-1}t,$ (\ref{ini}) is true as $(t-s)>2^{-1}t.$ So it is sufficient to consider the case
where $s\geq 2^{-1}t.$ From the form of $\GNRt$ in \eqref{GNR}, it is easy to see that our time-reversed process
has the following decomposition for fixed $t>0:$
\begin{eqnarray*}
& &\EP[ X^*_T(T-s)-X^*_T(T-t)|X_t=x]
= \EP\left[\int_s^t \GN f(X_r)dr\Big{|}X_t=x\right]\\ &+&
2\EP\left[\sum_{i,j}\int_s^t{\partial\over \partial x_j}\log h(0,r;x_0,X_r)g^{ij}(X_r){\partial\over \partial x_i}
f(X_r)dr\Big{|}X_t=x\right]\\
&\leq &(t-s)\| \GN f\|_\infty +C\EP\left[\int_s^t{|x_0-X_r|+c_2\over t^{c_2}}dr\Big{|}X_t=x\right]
\\
&\leq &(t-s)\| \GN f\|_\infty +C\int_s^ts^{-c_2}dr
\end{eqnarray*}
where we used Aronson's inequality in the last step, and
$C$ is a constant depending on the $C_1$-norm of $f$, $|x_0-x |,$ $c,$ $c_1$ and
$c_2.$  Dividing both sides of the last inequality by $(t-s)$ and noticing that
$s>2^{-1}t,$ we get the Lemma.   $\qed$
\end{pf*}
We also have
\begin{lem} \label{hz}
For all $f\in C^2_0(\Rd;\R)$ the following convergence
holds for almost every $y>0:$
\begin{equation*}
\begin{split}
&\lim_{t\searrow s}
\int_0^y\left\{(t-s)^{-{1\over 2}}p(s,t;y,z)
\left[\dstyle\int {h(0,z;x_0,\eta )h(z,y;\eta,\xi)\over h(0,y;x_0,\xi)}
f(\eta)d\eta -f(\xi )\right]\right\}dz
\\=&\ {A^*_yf(\xi )\over \sqrt{2 \pi }}.
\end{split}
\end{equation*}
Moreover, there is a constant $c_4$ such that
\begin{equation*}
\begin{split}
\int_0^y\left\{(t-s)^{-{1\over 2}}p(s,t;y,z)
\left[\dstyle\int {h(0,z;x_0,\eta )h(z,y;\eta,\xi)\over h(0,y;x_0,\xi)}
f(\eta)d\eta -f(\xi )\right]\right\}dz
\leq c_4y^{-c_2}
\end{split}
\end{equation*}
\end{lem}
\begin{pf*}{Proof}
By the reflection principle, the transition density of the reflecting BM $|B(s)|$ is
\begin{equation}
p(s,t; y,z)=\frac{1}{\sqrt{2\pi (t-s)}}\left[\exp\left\{{|y-z|^2\over 2(t-s)}\right\}
+\exp\left\{{|y+z|^2\over 2(t-s)}\right\}\right]
\end{equation}
By Lemma \ref{hhh},
\begin{equation*}
\begin{split}
&\int_0^y (t-s)^{-{1\over 2}}|y-z|p(s,t;y,z)
\frac{\left[\dstyle\int {h(0,z;x_0,\eta )h(z,y;\eta ,\xi )\over h(0,y;x_0,\xi )}
f(\eta)d\eta -f(\xi )\right]}{|y-z|}dz\\
\leq &\ C\int_0^y (t-s)^{-{1\over 2}}|y-z|p(s,t;y,z)  y^{-c_2}dz\\
\leq &\ C\int_0^y (t-s)^{-1}|y-z|\exp\left\{ -{|y-z|^2\over 2(t-s)}\right\}y^{-c_2}dz\\
=&\ C\int_0^{y/\sqrt{t-s}}z\exp\left\{ -{z\over 2}\right\}y^{-c_2}dz\\
\leq &\ Cy^{-c_2},
\end{split}
\end{equation*}
where $C$ is a generic constant that may vary from line to line.  Now, we may write
for $z<y$,
\begin{eqnarray*}
& &(y-z)^{-1}\left[\int {h(0,z;x_0,\eta )h(z,y;\eta ,\xi )\over h(0,y;x_0,\xi )}f(\eta )  d\eta -f(\xi )\right]
=\GNRy f(\xi )+o(y-z),
\end{eqnarray*}
where $o(y-z)\to 0$ (as $(y-z)\to 0$) and $o(y-z)\leq Cy^{-c_2}.$ On the other hand,
\begin{eqnarray*}
& &\lim_{t\searrow s}\int_0^y (t-s)^{-{1\over 2}}|y-z|p(s,t;y,z) [\GNRy f(\xi )+o(y-z)]dz\\
&=&\lim_{t\searrow s}\int_0^y
{|y-z|\over  \sqrt{2\pi }(t-s)}\exp\left\{-{|y-z|^2\over 2(t-s)}\right\}
[\GNRy f(\xi)+o(y-z)]dz\\
&=&\lim_{t\searrow 0}\int_0^{y/\sqrt{t-s}}
{|z|\over  \sqrt{2\pi }}\exp\left\{-{|z|^2\over 2}\right\}
[\GNRy f(\xi )+o(z\sqrt{t-s})]dz\\
&=&\GNRy f(\xi )\int_0^\infty
{|z|\over \sqrt{2\pi }}\exp\left\{-{|z|^2\over 2}\right\}dz\\
&=&{\GNRy f(\xi )\over \sqrt{2\pi }}.
\end{eqnarray*}
Thus we get the Lemma.  $\qed$
\end{pf*}
Similarly we have
\begin{lem} \label{hz2}
For all $f\in C^2_0(\Rd;\R)$ the following convergence
holds for almost every $y>0:$
\[ \lim_{t\searrow s}
\int_y^\infty (t-s)^{-{1\over 2}}p(s,t;y,z)
\left[\int h(y,z;\xi ,\eta )f(\eta ) d\eta -f(\xi )\right]dz={\GN f(\xi )\over
\sqrt{2 \pi }}   .\]
\end{lem}
\begin{pf*}{Proof of \thmref{HDGNthm}}
Now, we easily have
\[ \P(X(|B(s)|)\in d\xi )=\left[\int_0^\infty p(0,s;0,y)h(0,y;x_0,\xi )dy\right]d\xi .\]
And for $t>s,$ we see that

\begin{eqnarray*}
& & \P(X(|B(s)|)\in d\xi ,\ |B(t)|\geq |B(s)|, \  X(|B(t)|)\in d\eta )\\
&=& \left[\int_0^\infty \int_y^\infty p(0,s;0,y)p(s,t;y,z)
h(0,y;x_0,\xi )h(y,z;\xi ,\eta )dzdy\right]d\xi d\eta ,
\end{eqnarray*}
and
\begin{eqnarray*}
& & \P(X(|B(s)|)\in d\xi ,\ |B(t)|< |B(s)|, \   X(|B(t)|)\in d\eta )\\
&=& \left[\int_0^\infty \int_0^y p(0,s;0,y)p(s,t;y,z)
h(0,y;x_0,\xi )\P[X(z)\in d\eta | X(y)\in d\xi ]dzdy\right]d\xi  \\
&=& \left[\int_0^\infty \int_0^y p(0,s;0,y)p(s,t;y,z)
h(0,y;x_0,\xi ){h(0,z;x_0,\eta )h(z,y;\eta ,\xi )\over h(0,y;x_0,\xi )}dzdy\right]d\xi d\eta \\
&=& \left[\int_0^\infty \int_0^y p(0,s;0,y)p(s,t;y,z)
h(0,z;x_0,\eta )h(z,y;\eta ,\xi ) dzdy\right]d\xi d\eta.
\end{eqnarray*}
Thus,
\begin{eqnarray*}
& &\EP[f[X(|B_t|)]\ | \ X(|B_s|)=\xi]=\left\{ \int_0^\infty p(0,s;0,y)h(0,y;x_0,\xi )dy\right\}^{-1}
\\& &\left\{ \int \left[\int_0^\infty \int_y^\infty p(0,s;0,y)p(s,t;y,z)
h(0,y;x_0,\xi )h(y,z;\xi ,\eta )dzdy\right]f(\eta )  d\eta\right. \\
& &+\left.\int \left[\int_0^\infty \int_0^y p(0,s;0,y)p(s,t;y,z)
h(0,z;x_0,\eta )h(z,y;\eta ,\xi )dzdy\right]f(\eta )  d\eta \right\}
\end{eqnarray*}
\nopagebreak
and so to compute
$$\lim_{t\searrow s}\ (t-s)^{-{1\over 2}}\{ \EP\left[f(X(|B_t|))|X(|B_s|)\right]-f(X(|B_s|))
\},$$
we observe that
\begin{equation*}
\begin{split}
&\lim_{t\searrow s}\ {(t-s)^{-{1\over 2}}}\{\EP[ f[X(|B_t|)]\ | \ X(|B_s|)=\xi]-f(\xi)\}\\
=&\lim_{t\searrow s}\ (t-s)^{-{1\over 2}}\left\{\int_0^\infty p(0,s;0,y)h(0,y;x_0,\xi)dy
\right\}^{-1}\\
&\left\{\int \left[\int_0^\infty \int_y^\infty p(0,s;0,y)p(s,t;y,z)
h(0,y;x_0,\xi )h(y,z;\xi ,\eta )dzdy\right]f(\eta)d\eta\right.\\
&+\left.\int\left[\int_0^\infty \int_0^y p(0,s;0,y)p(s,t;y,z)
h(0,z;x_0,\eta )h(z,y;\eta ,\xi )dzdy\right]f(\eta)d\eta-f(\xi )\right\}\\
=&\lim_{t\searrow s}\ (t-s)^{-{1\over 2}}\left\{\int_0^\infty p(0,s;0,y)
h(0,y;x_0,\xi )dy\right\}^{-1}\\
&\left\{\int \left[\int_0^\infty \int_y^\infty p(0,s;0,y)p(s,t;y,z)
h(0,y;x_0,\xi )h(y,z;\xi ,\eta )dzdy\right]f(\eta)d\eta\right.\\
&+\int\left[\int_0^\infty \int_0^y p(0,s;0,y)p(s,t;y,z)
h(0,z;x_0,\eta)h(z,y;\eta,\xi)dzdy\right]f(\eta)d\eta\\
 &\left.-f(\xi)\int_0^\infty p(0,s;0,y)h(0,y;x_0,\xi )dy\right\}\\
=&\lim_{t\searrow s}\ (t-s)^{-{1\over 2}}
\left\{\int_0^\infty p(0,s;0,y)h(0,y;x_0,\xi)dy\right\}^{-1}\\
&\left\{\int\left[\int_0^\infty \int_y^\infty p(0,s;0,y)p(s,t;y,z)
h(0,y;x_0,\xi )h(y,z;\xi ,\eta )dzdy\right]f(\eta)d\eta\right.\\
&+\int\left[\int_0^\infty \int_0^y p(0,s;0,y)p(s,t;y,z)
h(0,z;x_0,\eta )h(z,y;\eta ,\xi )dzdy\right]f(\eta)d\eta\\
&\left.-f(\xi)\int_0^\infty \int_0^\infty p(0,s;0,y)h(0,y;x_0,\xi)p(s,t;y,z)dydz\right\}\\
\end{split}
\end{equation*}
\begin{equation}\label{hz3}
\begin{split}
=&\lim_{t\searrow s}\left\{\int_0^\infty p(0,s;0,y)h(0,y;x_0,\xi )dy\right\}^{-1}\\
&\left\{\int_0^\infty \int_y^\infty (t-s)^{-{1\over 2}}|y-z| p(0,s;0,y)p(s,t;y,z)
h(0,y;x_0,\xi )\right.\\
&|y-z|^{-1}\left[\int h(y,z;\xi,\eta)f(\eta)d\eta-f(\xi)\right]dzdy\\
&+\int_0^\infty \int_0^y (t-s)^{-{1\over 2}}|y-z|p(0,s;0,y)p(s,t;y,z)
h(0,y;x_0,\xi)\\
&\left.|y-z|^{-1}\left[\int {h(0,z;x_0,\eta)h(z,y;\eta,\xi)\over h(0,y;x_0,\xi)}f(\eta)d\eta-f(\xi)
\right]dzdy\right\}.
\end{split}
\end{equation}
It is easy to see by Lemma \ref{hz2} that
\begin{eqnarray*}
& &\lim_{t\searrow s} \left\{ \int_0^\infty p(0,s;0,y)h(0,y;x_0,\xi )dy\right\}^{-1}\\
& &   \left\{\int_0^\infty \int_y^\infty (t-s)^{-{1\over 2}}|y-z| p(0,s;0,y)p(s,t;y,z)
h(0,y;x_0,\xi )\right.\\& &\left.|y-z|^{-1}[\int h(y,z;\xi ,\eta )f(\eta)d\eta -f(\xi)]dzdy\right\} \\
&=&{1\over \sqrt{2 \pi }}\GN f(\xi ).
\end{eqnarray*}
So let us consider the last term in (\ref{hz3}). From Aronson's inequality (\ref{ar}) and
Lemma \ref{hz},
when $|x_0-\xi |>0,$
\begin{eqnarray*}
& &\int_0^y (t-s)^{-{1\over 2}}|y-z|p(s,t;y,z)
h(0,y;x_0,\xi )\nonumber\\
& &|y-z|^{-1}\left[\int {h(0,z;x_0,\eta )h(z,y;\eta ,\xi )\over h(0,y;x_0,\xi )}f(\eta )  d\eta -f(\xi )\right]dz\\
&\leq &c_4h(0,y;x_0,\xi )y^{-c_2}
\end{eqnarray*}
is bounded in $(t-s,y)$ for fixed $\xi$, and we may pass to the limit through
the integral over $\pR$.
Thus, the following half-derivative exists for every $s>0$ and is given by:
\begin{eqnarray*}
& &\lim_{t\searrow s}\ (t-s)^{-{1\over 2}}\{ \EP\left[f(X(|B_t|))|X(|B_s|)\right] -f(X(|B_s|)) \} \nonumber\\
& &={1\over \sqrt{2 \pi }}\left[\GN f(X(|B_s|) )+
{\frac{\dstyle\int_0^\infty p(0,s;0,y)h(0,y;x_0, X(|B_s|))\GNRy f(X(|B_s|) )dy}
{\dstyle\int_0^\infty p(0,s;0,y)h(0,y;x_0, X(|B_s|))dy}}\right],
\end{eqnarray*}
proving \thmref{HDGNthm}. $\qed$
\end{pf*}
\pagebreak
\appendix\section*{Appendix }
\setcounter{section}{1}
\setcounter{equation}{0}
We now rigorize and prove our claim in statement (c), in the introduction of this
paper, that $\EBTP$ is the $k\to\infty$ limit of $\kEBTP$.  This is accomplished by
showing weak convergence of the process $\left\{\kEBTP;0\le t<\infty\right \}$ to
$\left\{\EBTP ;0\le t<\infty\right\}$.
Without losing generality, we may assume that, for each $p>0$, there are positive constants $c_{1,p}$, $c_{2,p}$, and  $c_{3,p}$ such that
\begin{equation}\label{poth}
\begin{split}
 \P\left[\sup_{a\leq s\leq t\leq a+b}|X^{x,1}(t)-X^{x,1}(s)|^p>c_{1,p}b^{c_{2,p}p}\right]\leq \exp\left\{ -{c_{3,p}\over b}\right\}\quad\forall\ a,b\geq0.
\end{split}
\end{equation}
Clearly, \eqref{poth} is true when $\Xx$ is a Brownian motion, which is $\alpha$-H\"older continuous for any $\alpha<1/2$.   For a general $\Xx$, we see that
the  martingale part of the diffusion process $\Xx$ is of $\alpha$-H\"older continuous for any $\alpha<1/2$, and the non-martingale part is differentiable,
so it is even smoother, so \eqref{poth} is true here as well.   Now, note that the paths which do not satisfy
$$\sup_{a\leq s\leq t\leq a+b} |X^{x,1}(t)-X^{x,1}(s)|^p\leq c_{1,p}b^{c_{2,p}p}$$
have exponentially small probability, so they can be thrown away when $t-s$ is small.
\begin{thm}\label{th0}
There is a positive constant $c$ such that for each $p>0,$
there is a positive constant $C(p)$ satisfying
\begin{equation}\label{haha}
\EP\left|\kEBT(s)-\kEBTP\right|^p\leq C(p)|s-t|^{cp};\quad\forall\ 0\leq s\leq t<\infty, \forall k\in\N,
\end{equation}
and this is enough to conclude that there is a subsequence of $\left\{\kEBT\right\}$ converging weakly to $\EBT$, as $k\to\infty$.
\end{thm}
\begin{pf*}{Proof}
 Let $A_{i,s}\eqdef\left[ \kEBT(s)=X^{x,i}(|B(s)|)\right]$, for $1\le i\le k$ and
$0\le s<\infty$.  We then have
\begin{equation*}
\begin{split}
\EP\left|\kEBT(s)-\kEBTP\right|^p&=\sum_{i,j=1}^k\EP\left\{ 1_{A_{i,s}}1_{A_{j,t}}
\left|\kEBT(s)-\kEBTP\right|^p\right\}\\
&=\sum _{i,j=1\atop i\ne j}^{k}
\EP\left\{1_{A_{i,s}}1_{A_{j,t}}\left|X^{x,i}(|B(s)|)-X^{x,j}(|B(t)|)\right|^p\right\}\\
&+\sum_{i=1}^k\EP\left\{1_{A_{i,s}}1_{A_{i,t}}\left|X^{x,i}(|B(s)|)-X^{x,i}(|B(t)|)\right|^p\right\}\\
&=k(k-1)\EP\left\{1_{A_{s,1}}1_{A_{t,2}}
\left|X^{x,1}(|B(s)|)-X^{x,2}(|B(t)|)\right|^p\right\}\\
&+k\EP\{1_{A_{1,s}}1_{A_{1,t}}|X^{x,1}(|B(s)|)-X^{x,1}(|B(t)|)|^p\}\\
\end{split}
\end{equation*}
where the last equality follows from symmetry. From the definition of
$\kEBT(\cdot)$,
it is easy to see that the following inclusion of events is true when $i\neq j$:
\begin{equation*}
\left[\kEBT(s)=X^{x,i}(|B(s)|)\right]\cap\left[\kEBTP=X^{x,j}(|B(t)|)\right]
\subset
\left[\inf_{s\leq u\leq t}|B(u)| =0\right]\eqdef S_{s,t}.
\end{equation*}
Thus, by symmetry,
\begin{equation}\label{zhu}
\begin{split}
&\EP\left|\kEBT(s)-\kEBTP\right|^p\leq \EP\left\{1_{S_{s,t}}
\left|X^{x,1}(|B(s)|)-X^{x,2}(|B(t)|)\right|^p\right\}\\
&+k\EP\left\{ 1_{A_{1,s}}1_{A_{1,t}}\left|X^{x,1}(|B(s)|)-X^{x,1}(|B(t)|)\right|^p\right\}\\
&\leq C_p\EP\left\{1_{S_{s,t}}\left[\left|X^{x,1}(|B(s)|)-x\right|^p+\left|x-X^{x,2}(|B(t)|)\right|^p\right]\right\}\\
&+k\EP\left\{\left|X^{x,1}(|B(s)|)-X^{x,1}(|B(t)|)\right|^p\right\}
\end{split}
\end{equation}
As $x=X^{x,i}(0)$, then by (\ref{poth}) and the remarks following it and (\ref{zhu}), we obtain
\begin{equation}\label{zhuu}
\begin{split}
&\EP\left|\kEBT(s)-\kEBTP\right|^p\\
&\leq Cc_{1,p}\EP\left\{1_{S_{s,t}}\left[|B(s)|^{c_{2,p}p}
+|B(t)|^{c_{2,p}p}\right]\right\}+c_{1,p}\EP\left\{||B(s)|-|B(t)||^{c_{2,p}p}\right\}\\
&\leq Cc_{1,p}\EP\left\{ 1_{S_{s,t}}[|B(s)|^{c_{2,p}p}
+|B(t)|^{c_{2,p}p}]\right\}+c_{4,p}\EP\left\{|t-s|^{c_{5,p}p}\right\}
\end{split}
\end{equation}
where $C$ is a generic constant whose value may vary from line to line and
$c_{4,p}$ and $c_{5,p}$ are new constants obtained by
the well-known property of Brownian motion: there is a constant
${\overline C}_p$ such that
\begin{equation}\label{zhuuu}
\EP\left\{ \sup_{s_0\leq s\leq t\leq t_0}\left\{|B(t)-B(s)|^p\right\}\right\}
\leq {\overline C}_p|t_0-s_0|^{p\over 2};\qquad \forall\ 0\leq s\leq t<\infty.
\end{equation}
On the other hand, it is easy to see that
\begin{equation}\label{zh11}
\begin{split}
1_{S_{s,t}}\left\{|B(s)|^{c_{2,p}p}+|B(t)|^{c_{2,p}p}\right\}
\leq2\sup_{s\leq u\leq v\leq t}\{ |B(v)-B(u)|^p\}
\end{split}
\end{equation}
Thus, (\ref{haha}) can be easily deduced from (\ref{zhuu}), (\ref{zh11}), and (\ref{zhuuu}).
\vskip .2in
It is well known (see, e.g. \cite{El} and \cite{Ku}) that Kolmogorov's criterion implies that the sequence
of processes $\left\{\kEBTP; 0\le t<\infty\right\}_k $ is tight in law under the uniform convergence topology.
It is easy to check that any limit of the convergent subsequence of $\left\{\kEBT\right\}$ gives the
law of $\EBT$.  Thus we proved statement (c) in Section 0.  
\end{pf*}
 %%%%%%%%%%%%

    \section*{Acknowledgments.}
     The authors are deeply indebted to Chris Burdzy for linking them up in this project.
They also like to thank Pat Fitzsimmons for bringing \cite{Kim} to their attention.
The first author would like to thank Chris Burdzy, Davar Khoshnevisan, and Rich Bass for their
constant encouragements and Rick Durrett for fruitful discussions at the very early
stages of this project.  The authors would also like to thank an anonymous referee for a careful reading of our paper and for his useful comments
which improved it.
   %

          \iffalse
\begin{acknowledge}
The authors are deeply indebted to Chris Burdzy for linking them up in this project.
They also like to thank Pat Fitzsimmons for bringing \cite{Kim} to their attention.
The first author would like to thank Chris Burdzy and Davar Khoshnevisan for their
constant encouragements and Rick Durrett for fruitful discussions at the very early
stages of this project.  The authors would also like to thank an anonymous referee for a careful reading of our paper and for his useful comments
which improved it.
\end{acknowledge}
  \fi
  %

%%%%%%%%%%%%%%
\end{document}